\documentclass[11pt,reqno]{amsart}
\usepackage{amsmath,amsfonts,  amssymb,amsthm,amscd}
\makeindex

\usepackage{graphics}

\usepackage[latin1]{inputenc}
\usepackage{epsf}

\usepackage{enumerate}
\usepackage{indentfirst}
\usepackage{euscript}
\usepackage{graphicx}
\usepackage{amsmath}
\usepackage{amsfonts}
\usepackage{amssymb}
\makeindex

 \theoremstyle{plain}
\newtheorem{proposition}{Proposition}
\newtheorem{lemma}{Lemma}

\newtheorem{remark}{Remark}
\theoremstyle{definition}

\newtheorem{theorem}{Theorem}

 \newcommand{\ep}{\epsilon}

 \title[ Geometric Mean Curvature Lines ]{ Geometric Mean Curvature Lines  on Surfaces Immersed  in $\mathbb R^3$ }

 \author[R. Garcia]{Ronaldo Garcia}

\author[J. Sotomayor]{Jorge Sotomayor}


 \thanks{The first author was partially supported by     FUNAPE/UFG.
Both   authors are fellows of  CNPq.
 This work was done under the project PRONEX/FINEP/MCT - Conv. 76.97.1080.00 - Teoria Qualitativa das Equa\c c\~oes Diferenciais Ordin\'arias and 
CNPq - Grant 476886/2001-5.}

\begin{document}
 \maketitle

 {\small {\sc Abstract.-}
Here are studied  pairs of transversal  foliations
 with singularities, defined on  the  {\it Elliptic} region (where the Gaussian curvature  $\mathcal K$ is positive) of an oriented surface immersed in
$\mathbb R^3$. The   leaves of the foliations are the  {\it lines of  geometric mean
curvature,}  along which   the normal curvature is given by $\sqrt {\mathcal K}$, which is the  {\it geometric mean
curvature} of the  principal curvatures ${k_1}, k_2$ of the immersion.

The singularities of the foliations  are
 the {\it umbilic points} and {\it parabolic curves}, where ${k_1} = k_2$ and  ${\mathcal K} = 0$, respectively.

Here are determined the structurally  stable patterns of geometric mean curvature lines  near the
 {\it umbilic points},   {\it parabolic curves}  and {\it  geometric mean  curvature  cycles}, the periodic  leaves of the foliations. The genericity of these patterns is established.

This  provides  the three essential local ingredients to establish
 sufficient
conditions, likely to be also  necessary,  for {\it Geometric  Mean Curvature Structural Stability}. This study, outlined at the end of the paper, is   a natural   analog
and
complement  for the   {\it  Arithmetic Mean Curvature}  and {\it Asymptotic Structural Stability} of immersed surfaces studied previously by the authors \cite {a1, m, a2}.

\vskip .4cm

  {\sc R\'esum\'e.} -
Dans ce travail on \'etudie  les
  paires de
  feuilletages transverses avec  singularit\'es,
  d\'efinis dans la
  r\'egion  elliptique d'une surface orient\'ee
  plong\'ee  dans l'espace
  $\mathbb R^3$.
  Les feuilles sont  les lignes de  {\it courbure
  g\'eom\'etrique   moyenne,} selon lesquelles la
  courbure normale  est
  donn\'ee par la moyenne g\'eom\'etrique  $\sqrt {k_1
  k_2}$ des courbures
  principales  ${k_1}, k_2$.
  Les singularit\'es sont les points ombilics ( o\'u
  ${k_1} = k_2$)  et les
  courbes
  paraboliques (o\'u ${k_1}{k_2}  =0$ ).

  On d\'etermine  les conditions pour la
  stabilit\'e structurelle des
  feuilletages autour  des points ombilics, des
  courbes
  paraboliques et
  des cycles de courbure g\'eom\'etrique moyenne  (qui
  sont les feuilles
  compactes).
  La g\'en\'ericit\'e de ces conditions est
  \'etablie.

  Munis de  ces conditions    on
  \'etablit  les conditions suffisantes,  qui sont aussi
  vraisemblablement
   n\'ecessaires,  pour la stabilit\'e
  structurelle des
  feuilletages. Ce travail est  une continuation et une
g\'en\'eralisation
   naturelles de
  ceux  sur les lignes \`a {\it courbure  arithmétique
  moyenne}, selon lesquelles la courbure normale est
   donn\'ee par $(k_1+k_2)/2$ et sur les lignes  \`a {\it
  courbure nulle,} qui sont les   {\it courbes
  asymptotiques}. Voir les  articles   \cite{a1, m, a2}.


\section{Introduction}

In this  paper are  studied  the {\it geometric mean curvature
configurations} associated to  immersions of  oriented surfaces
into $\mathbb R^3$. They consist on the {\it umbilic points} and
{\it parabolic curves}, as singularities,   and of  the {\it lines
of geometric mean  curvature}  of the immersions, as the leaves of
the two transversal foliations in the configurations.  Along these
lines   the  normal curvature is  given by the geometric mean
curvature, which is the   square root of the  product of the
principal curvatures (i.e of the Gaussian Curvature).

The two transversal foliations, called here {\it geometric mean curvature foliations}, are well defined and regular only on the non-umbilic part of the elliptic region of the immersion, where the  Gaussian Curvature is positive. In fact, there they are the solution of smooth quadratic differential equations.  The set where the Gaussian Curvature vanishes,  the parabolic set, is generically a  regular curve which is the border of the elliptic region. The umbilic points are those at which the principal curvatures coincide,  generically are isolated and disjoint from the parabolic curve. See section 2 for precise definitions.

This study is  a natural development and extension of  previous results  about  the  Arithmetic Mean Curvature and Asymptotic
Configurations,
 dealing with the qualitative properties of the  lines along which the  normal  curvature is  the  arithmetic mean of the principal curvatures (i.e. is  the standard Mean Curvature) or is null.
This  has been considered  previously by the authors; see  \cite {a1, a2} and  \cite {m}.

The point of departure of this line of research, however, can be found in the classical works of Euler,  Monge, Dupin and Darboux, concerned with the  lines of principal curvature and umbilic points  of immersions. See \cite {Sp, St} for an initiation on  the basic facts on this subject; see  \cite {gs1, gs2} for a  discussion of the classical contributions and for their analysis  from the point of view of  structural stability  of differential equations  \cite {pm}.

This paper establishes  sufficient conditions, likely to be also necessary, for the  structural stability of {\it geometric mean curvature configurations}
  under small perturbations of the immersion. See section 7  for precise statements.

  This  extends to the  geometric mean curvature setting the  main theorems on structural stability   for the arithmetic mean curvature
configuration and for the asymptotic configurations, proved in \cite {a1, m, a2}.

Three  local ingredients  are essential  for this extension: the umbilic
points, endowed  with  their geometric  mean curvature separatrix structure, the geometric
mean curvature cycles, with the calculation of the derivative of the Poincar\'e return map,  through which is expressed the hyperbolicity condition and   the parabolic curve, together with the parabolic tangential  singularities
and associated separatrix structure.

The conclusions  of this paper, on the elliptic  region, are  complementary
to  results valid independently on the
hyperbolic region (on which the Gaussian curvature is negative), where the separatrix  structure near the parabolic curve and the asymptotic structural stability has been studied in \cite {a1, a2}.

The parallel  with  the conditions for principal, arithmetic  mean curvature and asymptotic  structural stability is remarkable. This can be attributed to the unifying role played by the notion of Structural Stability  of  Differential Equations and Dynamical Systems, coming to Geometry through the seminal work of Andronov and  Pontrjagin \cite {ap} and Peixoto \cite {mp}.

  The interest on lines of geometric mean curvature  goes back to the paper  of Occhipinti \cite {oc}.  The work of  Ogura \cite {og} regards these lines in terms   of his unifying notions {\it T-Systems} and {\it K-Systems} and makes a local analysis of the expressions of the fundamental quadratic forms in a chart whose coordinate curves  are   lines of geometric mean curvature. A  comparative study of these expressions with those corresponding to other lines  of geometric interest, such as  the {\it principal, asymptotic, arithmetic mean curvature}  and {\it characteristic lines}, was carried out by Ogura in the  context of {\it T-Systems} and {\it K-Systems}.  In  \cite{h}  the authors have studied the foliations by {\it characteristic lines}, called  harmonic mean curvature lines.

The authors are grateful to Prof. Erhard Heil for calling their attention to these papers, which seem to have remained  unquoted along so many years.

No global examples, or even local ones around singularities, of geometric mean curvature configurations seem to have been considered in the  literature on differential equations of classic differential geometry, in contrast with the situations for the   principal and asymptotic cases mentioned above. See also the work of  Anosov, for the global structure of the geodesic flow \cite {an},   and that  of Banchoff, Gaffney and McCrory \cite {bgm} for the parabolic and asymptotic lines.

\vskip 0.2cm

This paper is organized as follows:

Section 2 is devoted to the general study of the differential equations and general properties of Geometric Mean Curvature Lines. Here are given the precise definitions of the Geometric Mean Curvature Configuration  and of the  two transversal Geometric Mean Curvature Foliations with singularities into which it splits.  The definition of Geometric Mean Curvature Structural Stability focusing the preservation  of the qualitative properties of the foliations and the configuration under small perturbations of the immersion, will be given at the end of this section.

In Section 3 the equation of lines of geometric mean curvature is written in a Monge chart. The condition  for umbilic geometric mean curvature stability  is explicitly stated in terms of the coefficients of the third order jet of the    function which represents the immersion in a Monge chart. The local geometric mean curvature separatrix configurations at stable umbilics is established for $C^4$ immersions and resemble the three  Darbouxian patterns of principal and arithmetic mean curvature configurations \cite {da, gs1}.

In Section 4 the derivative of first return Poincar\'e map along
a geometric mean curvature cycle is established. It consists of an integral
expression
involving the curvature  functions
along the  cycle.

In Section 5 are studied the foliations by  lines of geometric
mean curvature  near the parabolic set   of an immersion, assumed
to be a regular curve. Only two generic patterns of the three
singular tangential patterns in common with the asymptotic
configurations, the {\it folded node} and the  {\it folded
saddle},  exist generically in the case; the {\it folded focus}
being  absent. See \cite{a1}.

Section 6 presents examples of Geometric Mean Curvature Configurations on  the Torus of revolution  and  the quadratic Ellipsoid, presenting non-trivial recurrences. This situation, impossible in principal configurations, has been established for arithmetic mean curvature configurations in \cite{m}.

In Section 7 the results presented in Sections 3, 4 and 5 are put together
to provide sufficient conditions for Geometric Mean Curvature Structural Stability. The genericity of these conditions is formulated at the end of this section, however its rather technical proof will be postponed to  another paper.

 \section{Differential Equation of Geometric Mean Curvature Lines}

Let $\alpha : {\mathbb M}^2\to \mathbb R^3$ be a $C^r,\;\; r\geq 4,$ immersion of
an oriented smooth surface ${\mathbb M}^2$ into $\mathbb R^3$. This means that $D\alpha$
is injective  at every point in ${\mathbb M}^2$.

 The space $\mathbb R^3$ is oriented by a  once for all fixed orientation
and endowed with the Euclidean inner product  $<,>$.

Let $N$  be a vector field orthonormal to $\alpha$.
Assume that $(u,v)$ is a positive chart of ${\mathbb M}^2$  and that $\{\alpha_u, \alpha_v, N\}$
is  a positive frame in $\mathbb R^3$.

In the  chart $(u,v)$, the {\it first fundamental form} of an immersion
$\alpha$ is
given by:

$I_\alpha= <D\alpha,D\alpha>= E du^2+2Fdudv+Gdv^2$, with

$E=<\alpha_u,\alpha_u>$, $F=<\alpha_u,\alpha_v>$, $G=<\alpha_v,\alpha_v>$

The  {\it second fundamental form} is given by:

 \centerline{$II_{\alpha}= <N , D^2\alpha> = e du^2 + 2f dudv + g dv^2$.}

The normal curvature  at a point $p$ in a tangent direction $t=[du:dv]$ is given by:

 $$ k_n= k_n(p) =\frac{ II_\alpha(t,t)}{I_\alpha(t,t)}.$$

The lines of geometric   mean curvature  of $\alpha$ are regular curves $\gamma$ on $\mathbb M^2$  having
normal  curvature   equal to the
geometric mean curvature of the immersion, i.e.,  $k_n=\sqrt{\mathcal K}$, where $ {\mathcal K}$ = ${\mathcal K}{_\alpha}$ is the Gaussian curvature of $\alpha$.

Therefore the   pertinent differential equation for these lines is given by:

$$\frac{edu^2+2fdudv+gdv^2}{Edu^2+2Fdudv+Gdv^2}=\sqrt{\frac{eg-f^2}{ EG-F^2}}=\sqrt{\mathcal K}$$

Or equivalently by

\begin{equation}\label{eq:mgc}
 [g-\sqrt{\mathcal K}G]dv^2+2[f-\sqrt{\mathcal K}F]dudv+  [e-\sqrt{\mathcal K}E]du^2=0.\end{equation}

This equation is defined only on the closure of the {\it Elliptic region}, ${\mathbb E}{\mathbb M^2}{_\alpha}$,   of $\alpha$, where ${\mathcal K} >0$. It is bivalued and $C^{r-2},\;\; r\geq 4,$ smooth on the complement of the umbilic,  ${\mathcal U}_\alpha$, and parabolic, ${\mathcal P}_\alpha$, sets of the immersion $\alpha$.  In fact, on ${\mathcal U}_\alpha$ , where the principal curvatures  coincide, the equation vanishes identically; on  ${\mathcal P}_\alpha$, it is univalued.

Also, the above equation is equivalent to the quartic differential equation, obtained from the above one by eliminating the square root.

\begin{equation}\label{eq:mgc4}A_{40} du^4+ A_{31}du^3dv + A_{22}du^2dv^2+A_{13} dudv^3+A_{04} dv^4 =0
\end{equation}

\noindent where,

$$\aligned A_{40}&= e^2 (E G-F^2)-E^2 (eg-f^2)\\
A_{31}&=4ef(EG-F^2)-4EF(eg-f^2)\\
A_{22}&= 6f^2 EG-6eg F^2\\
A_{13}&=4fg(EG- F^2)-4FG(eg-f^2) \\
A_{ 04}&= g^2 (E G-F^2)-G^2 (eg-f^2) \endaligned$$

\vskip 0.3cm

 The developments above allow us
to organize the lines of geometric mean curvature of immersions  into  the
{\it geometric mean curvature
configuration,}  as follows:

 Through every point $p\in {\mathbb E}{\mathbb M^2}{_\alpha}\backslash({\mathcal U}_\alpha \cup {\mathcal P}_\alpha)$,  pass two
geometric   mean curvature
lines of $\alpha$. Under the orientability
hypothesis imposed on $\mathbb M$, the geometric mean curvature
lines define two foliations:  $\mathbb G_{\alpha,1}$,
called the {\it minimal geometric mean curvature foliation}, along which the
geodesic torsion is
negative (i.e  $\tau_g=-\sqrt[4]{\mathcal K}\sqrt{ 2{\mathcal H}-2\sqrt{{\mathcal K}}}$ ),  and $\mathbb G_{\alpha,2}$, called the {\it maximal
geometric mean curvature foliations}, along which the geodesic torsion is
positive  (i.e  $\tau_g=\sqrt[4]{\mathcal K}\sqrt{ 2{\mathcal H}-2\sqrt{{\mathcal K}}}$  ).

By comparison with the arithmetic mean curvature directions, making angle $\pi/4$ with the minimal principal directions, the geometric ones are located  between them and the principal ones, making an angle $\theta$ such that
$tan\theta$ =$\pm\sqrt[4]{\frac{k_1}{k_2}}$,
 as follows from Euler's Formula. The particular expression for the geodesic torsion given above
results from the expression $\tau_g=(k_2-k_1)sin{\theta} cos{\theta}$ \cite{St}.  It is found in the work of Occhipinti \cite {oc}. See also Lemma 1 in Section 4 below.

The quadruple
$\mathbb G_\alpha=\{{\mathcal P}_\alpha,
{\mathcal U}_\alpha, \mathbb G_{\alpha,1},\mathbb G_{\alpha,2}\}$
is called the {\it geometric mean curvature configuration} of $\alpha$.
It splits into two foliations with singularities:
$\mathbb G^i_\alpha=\{{\mathcal P}_\alpha,
{\mathcal U}_\alpha, \mathbb G_{\alpha,i}\}, \; i=1, 2$.

Let $\mathbb M^2$ be also   compact. Denote by ${\mathcal M}^{r,s}({\mathbb M^2)}$ be the space of $ C^r $ immersions of $\mathbb M^2$ into the Euclidean space $\mathbb R^3$, endowed with the $C^s$ topology.

 An immersion $\alpha$ is said $C^s$-{\it local  geometric mean curvature structurally
stable} at a compact set  $K\subset \mathbb M^2$
if for any sequence of immersions $\alpha_n$ converging to $\alpha$ in ${\mathcal M}^{r,s}({\mathbb M^2)}$, in a compact neighborhood $V_K$ of $K$, there is a sequence of compact subsets $K_n$ and a sequence of homeomorphisms mapping $K$
to $K_n$, converging to the identity of $\mathbb M^2$, such that on $V_K$ it  maps umbilic,   parabolic  curves and  arcs of the geometric mean curvature foliations
$\mathbb G_{\alpha,i}$ to those  of $\mathbb G_{\alpha_n,i}$
for $i = 1,\;2$.

An immersion $\alpha$ is said $C^s$-{\it geometric mean curvature structurally stable} if the compact $K$ above is the closure of  ${\mathbb E}{\mathbb M^2}{_\alpha}$.

Analogously, $\alpha$ is said {\it i-} $C^s$-{\it geometric mean curvature structurally stable}
if only the preservation of elements of  {\it i-th, i=1,2} foliation with singularities is required.

A general study of the structural stability of  quadratic
differential equations (not necessarily derived from normal
curvature properties) has been carried out by Gu{\'\i}\~nez \cite
{gn}. See also the work of Bruce and Fidal  \cite {bf} for the analysis of umbilics for general quadratic differential equations.

\section{ Geometric mean curvature  lines near umbilic points }

 Let  $0$ be an umbilic point of a $ C^r ,\; r\geq  4,$ immersion  $\alpha$ parametrized in a Monge chart $(x,y)$ by $\alpha(x,y)=(x,y,h(x,y))$, where

\begin{equation} \label{eq:1} h(x,y)=\frac k2(x^2+y^2)+ \frac a6 x^3+\frac b2 xy^2+\frac c6 y^3+O(4)\end{equation}

This reduced form is obtained by means of a rotation of the $x,y$-axes. See \cite {gs1, gs2}.

According to Darboux \cite{da, gs1},  the differential equation of principal curvature lines  is given by:

\begin{equation} -[by+P_1]dy^2+[(b-a)x+cy+P_2 ]dxdy+[by+P_3 ]dx^2=0, \end{equation}

\noindent Here the  $P_i=P_i(x,y)$  are functions of the form $P_i(x,y)=O(x^2+y^2)$.

 As an starting point,   recall the behavior of principal lines near
Darbouxian umbilics
 in the following proposition.

\begin{proposition}{\rm \cite {gs1, gs2}} \label{prop:1} Assume  the notation  established in \ref{eq:1}.
Suppose that the transversality condition $T:  b(b-a)\ne 0$  holds
 and
consider the following situations:
\begin{itemize}
\item[$D_1$)] $\;\;\Delta_{P}>0$

\item[$D_2$)] $\;\;   \Delta_{P} <0$ and $\dfrac ab >1$

\item[$D_3$)]   $\;\;    \dfrac ab  <1  $
\end{itemize}
 Here $\Delta_{P} =4b(a-2b)^3-c^2(a-2b)^2$

 Then each principal foliation  has  in a neighborhood of $0$,  one
hyperbolic sector in the $D_1$ case, one
parabolic and one hyperbolic sector in $D_2$ case and three hyperbolic sectors in the case $D_3$.
  These  points  are called principal curvature  Darbouxian umbilics.
\end{proposition}

\begin{proposition}\label{prop:2}  Assume  the notation  established in \ref{eq:1}.
Suppose that the transversality condition $T_g :  kb(b-a)\ne 0$  holds
 and
consider the following situations:
\begin{itemize}
\item[$G_1$)]   $\;\;\Delta_{G}>0$

\item[$G_2$)]$\;\;\Delta_{G}<0 $ and $ \;\;  \dfrac ab  >1$

\item[$G_3$)]    $\;\; \dfrac ab < 1. $
\end{itemize}
Here  $ \Delta_{G} =4c^2(2a-b)^2-[3c^2+(a-5b)^2][3(a-5b)(a-b)+c^2].  $

  {Then each geometric mean curvature foliation  has  in a neighborhood of $0$,  one
hyperbolic sector in the $G_1$ case, one
parabolic and one hyperbolic sector in $G_2$ case and three hyperbolic sectors in the case $G_3$. These umbilic points are called geometric mean curvature Darbouxian umbilics.}

 The  geometric mean curvature foliations $\mathbb G_{\alpha,i}$ near an umbilic point of type $G_k$ has  a  local behavior  as shown in Figure 1. The separatrices of these singularities are called umbilic separatrices.
\end{proposition}

\begin{figure}[htbp] \label{fig:dgeo}
\begin{center}
\includegraphics[angle=0, width=11cm]{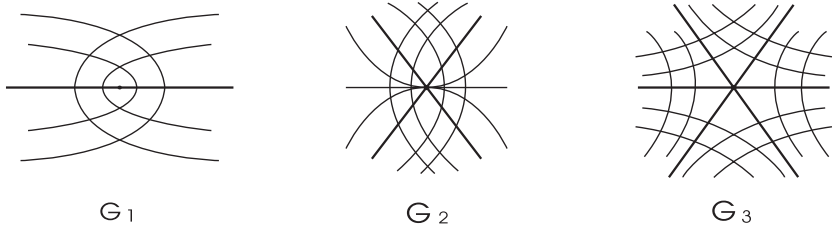}
       \caption {  Geometric mean curvature lines near the umbilic points $G_i$ and their separatrices  }
   \end{center}
 \end{figure}

\begin{proof}
Near $0$ the functions $\sqrt{\mathcal K}$ and $\mathcal H$ have the following Taylor expansions. Assume here that $k>0$, which can be achieved by means of an exchange in orientation.
So it follows that,

$$\sqrt{\mathcal K}=k+ \frac 12(a+b)x+\frac 12cy+O_1(2),\hskip .5cm
{\mathcal H}=k+ \frac 12(a+b)x+\frac 12 cy+O_2(2).$$

Therefore the  differential equation of the geometric mean curvature lines $$ [g-\sqrt{\mathcal K}G]dv^2+2[f-\sqrt{\mathcal K}F]dudv+  [e-\sqrt{\mathcal K}E]du^2=0 $$\noindent  is given by:

\begin{equation} \label{eq:mg}[(b-a)x+cy+M_1]dy^2+ [4by+M_2]dxdy- [(b-a)x+cy+M_3]dx^2 =0\end{equation}

 \noindent where $M_i$,  $i=1,2,3$,  represent  functions of
order $O((x^2+y^2)).$

At the level of first jet the differential equation \ref{eq:mg} is the same as that of the arithmetic mean curvature lines given by
$$ [g- {\mathcal H}G]dv^2+2[f- {\mathcal H}F]dudv+  [e- {\mathcal H}E]du^2=0.$$

The condition $\Delta_G$ coincides with  the $\Delta_H$ condition established to characterizes the arithmetic mean curvature Darbouxian umbilics studied in detail  in \cite{m} reducing the analysis of that of hyperbolic saddles and nodes whose phase portrait is determined   only by  the first jet.
\end{proof}

\begin{remark}\label{prop:mgumbilico}In the plane $b=1$ the bifurcation diagram of the  umbilic points of types $G_i$ for the geometric mean curvature configuration and of types $D_i$ for the principal configuration is as shown in
 Figure 2  below.

\end{remark}

\begin{figure}[htbp]
\begin{center}
\includegraphics[angle=0,  width=7cm]{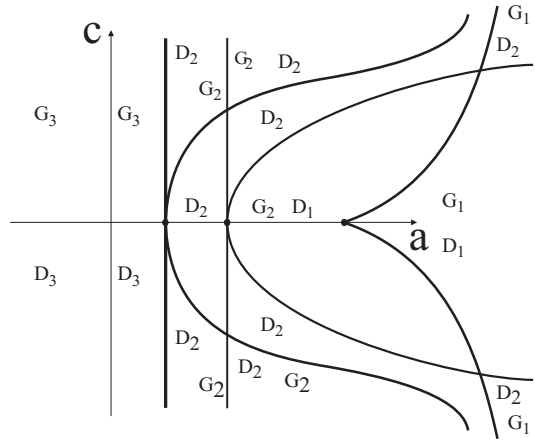}
      \caption {Bifurcation diagram of points $D_i$ and $G_i$}
   \end{center}
\end{figure}

\begin{theorem}\label{th:31}
An immersion $\alpha\in {\mathcal M}^{r,s}({\mathbb M^2)}$, $ r\geq   4$, is $C^3-$local
geometric mean curvature structurally  stable at ${\mathcal U}_\alpha$ if and only if
every $p\in {\mathcal U}_\alpha$ is one of the types $G_k$, $k=1,2,3$ of proposition \ref{prop:2}
\end{theorem}

\begin{proof} Clearly proposition \ref{prop:2} shows that the conditions $G_i$, $i=1,2,3$ and $T_g:kb(b-a)\ne 0 $ included  imply
the $C^3-$local geometric mean curvature structural stability. This involves the construction of  the homeomorphism (by means of
canonical regions) mapping simultaneously  minimal and maximal geometric mean curvature lines  around the umbilic points of
$\alpha$ onto those of a $C^4$ slightly perturbed immersion.

We will discuss the necessity of the  condition $T_g:k(b-a)b\ne
0$ and of the conditions $G_i$, $i=1, \,2,\, 3$.  The first one  follows from
its identification with a transversality condition that guarantees the persistent isolatedness  of the umbilic
points of $\alpha$ and its separation from the parabolic set, as well as the persistent regularity of  the Lie-Cartan
  surface
${\mathcal G}$.  Failure of $T_g$
condition has the following implications:
\begin{itemize}
\item[a)] $b(b-a)=0$; in this case the elimination or splitting of the umbilic point can be achieved by small perturbations.

\item[b)] $k=0$ and $b(b-a)\ne 0$; in this case a small perturbation separates the umbilic point from the parabolic set.
\end{itemize}
 The necessity of $G_i$ follows from its dynamic identification with the hyperbolicity  of the
equilibria along the projective line of the vector field   ${\mathcal G}$. Failure of this condition would  make possible
to  change the number of geometric mean curvature umbilic separatrices at the umbilic point by means a small  perturbation of
the immersion.
\end{proof}

\section{Periodic Geometric Mean Curvature Lines}

Let $\alpha:\mathbb M^2\to \mathbb R^3$ be an immersion of a compact and
 oriented surface and consider the foliations $\mathbb G_{\alpha, i}$, $i=1,\;2$, given by
 the {\it geometric mean curvature lines}.

  In terms of geometric invariants, here is established an integral
expression  for  the first derivative of the  return map of a
periodic   geometric mean curvature line, called {\it geometric mean curvature cycle}.
Recall that the return map associated  to a cycle
is a local diffeomorphism with a fixed point, defined on a cross section  normal to the cycle by following  the
integral curves through this section until they meet again the section. This map is called holonomy in   Foliation
Theory and Poincar\'e Map in Dynamical Systems, \cite{pm}.

A geometric mean curvature cycle is called {\it hyperbolic} if the first derivative of the return map at the fixed  point is
different from one.

{ The geometric mean curvature foliations $\mathbb G_{\alpha,i}$ has no geometric mean curvature cycles such that  the return map reverses the orientation.}
 Initially, the integral expression for the derivative of the return
map is
obtained in class $ C^6$; see Lemma \ref{lm:42} and Proposition \ref{prop:41}. Later on,
in Remark
\ref{rm:41} it is shown how to extend it to class $C^3$.

The characterization of hyperbolicity of geometric mean curvature cycles in terms
of local structural stability is given in  Theorem  \ref{th:41} of
this section.

\begin{lemma}\label{lm:41} Let $c:I\to {\mathbb M}^2$ be a geometric mean curvature
line parametrized by arc length. Then the Darboux frame is given by:

$$\aligned
 T^\prime &= k_g N\wedge T+ \sqrt{\mathcal K} N\\
(N\wedge T)^\prime &= -k_g T+ \tau_g N\\
N^\prime &= -\sqrt{\mathcal K} T-\tau_g N\wedge T\endaligned\eqno
$$
 \noindent where $\tau_g= \pm \sqrt{2{\mathcal H} -2\sqrt{\mathcal K}} \sqrt[4]{{\mathcal K}} $. The sign of $\tau_g$
is positive (resp. negative) if $c$ is maximal (resp. minimal) geometric mean curvature line.
\end{lemma}

\begin{proof} The normal curvature $k_n$ of the curve $c$ is by the definition
 the  geometric mean curvature $\sqrt{\mathcal K}$.
From the Euler equation $k_n=k_1\cos^2\theta+k_2\sin^2\theta=\sqrt{\mathcal K}$,
get $\tan\theta=\pm\sqrt{\frac{\sqrt{\mathcal K}-k_1}{k_2-\sqrt{\mathcal K}}}$.
    Therefore, by direct calculation, the geodesic torsion is given by
$\tau_g=(k_2-k_1)\sin\theta\cos\theta  =  \pm \sqrt{2{\mathcal H} -2\sqrt{\mathcal K}} \sqrt[4]{{\mathcal K}}$.
\end{proof}

 \begin{remark}
The expression for the geodesic curvature $k_g$ will not be needed explicitly in this work. However, it can be given in terms of the principal curvatures and their derivatives using a formula due to Liouville \cite{St}.
\end{remark}

\begin{lemma}\label{lm:42} Let $\alpha:\mathbb M\to \mathbb R^3$ be an immersion of class $ C^r $, $ r\geq   6$, and $c$ be a mean
curvature cycle of $\alpha$,
 parametrized by arc length $s$ and  of length $L$. Then the expression,

$$\alpha(s,v)=c(s)+v(N\wedge T)(s)+[\frac{2{\mathcal H}(s)-\sqrt{{\mathcal K}(s)}}{2}v^2+
 \frac{A(s)}{6}v^3+v^3B(s,v)]N(s)\eqno$$

 \noindent where $B(s,0)=0$, defines    a local chart $(s,v)$ of class $C ^{r-5}$ in a neighborhood of $c$.

\end{lemma}

\begin{proof}  The curve $c$ is of class $C ^{r-1}$ and the map
$\alpha(s,v,w)=c(s)+ v(N\wedge T)(s)+wN(s)$ is of class $C ^{r-2}$ and is a local diffeomorphism  in a neighborhood of the axis $s$. In fact $[\alpha_s,\alpha_v,\alpha_w](s,0,0)=1$. Therefore there is a function $W(s,v)$ of class $C ^{r-2}$ such that $\alpha(s,v,W(s,v))$ is a parametrization of a tubular neighborhood of $\alpha\circ c$. Now for each $s$, $W(s,v)$ is just a parametrization of
the curve of intersection between $\alpha(\mathbb M)$ and  the normal plane generated by $\{(N\wedge T)(s), N(s)\}$. This curve of intersection  is tangent to $(N\wedge T)(s)$ at $v=0$ and notice that  $k_n(N\wedge T)(s)=2{\mathcal H}(s)-\sqrt{{\mathcal K}(s)}$. Therefore,

\begin{equation} \aligned  \alpha(s,v,W(s,v))=& c(s)+ v(N\wedge T)(s)\\
+&[\frac{2{\mathcal H}(s)-\sqrt{{\mathcal K}(s)}}{2}v^2+
 \frac{A(s)}6v^3+ v^3B(s,v)]N(s), \endaligned \end{equation}

 \noindent where $A$ is of class $C ^{r-5}$ and $B(s,0)=0$.
\end{proof}

We now  compute the coefficients of the first and second fundamental forms in  the chart $(s,v)$ constructed above, to be used in  proposition \ref{prop:41}.

$$\aligned N(s,v)=&\frac{\alpha_s\wedge\alpha_v}{\mid\alpha_s\wedge\alpha_v\mid}  = [-\tau_g(s)v+O(2)]T(s)\\
-&[(2{\mathcal H}(s)-2\sqrt{{\mathcal K}(s)})v+O(2)](N\wedge T)(s)
+[1+ O(2)]N(s).\endaligned $$

 \noindent Therefore it follows that  $E=<\alpha_s,\alpha_s>$, $F=<\alpha_s,\alpha_v>$, $G=<\alpha_v,
\alpha_v>$,
$e=<N,\alpha_{ss}>,$
 $\;\; f=<N,\alpha_{sv}>\;\;$ and $\; g=<N,\alpha_{vv}>$ are given by

\begin{equation}\label{eq:1f2f}
\aligned
E(s,v) &= 1-2k_g(s)v+h.o.t\\
F(s,v) &=  0+ 0.v+h.o.t\\
G(s,v) &= 1 +0.v+h.o.t\\
e(s,v)&=\sqrt{{\mathcal K} (s)}+v[\tau_g^\prime(s)-2k_g(s) {\mathcal H}(s) ]+ h.o.t\\
f(s,v) &= \tau_g(s)+ \{[2{\mathcal H}(s)-\sqrt{{\mathcal K}(s)}]^\prime +k_g(s)\tau_g(s)\}v+ h.o.t\\
g(s,v) &=2{\mathcal H}(s)-\sqrt{{\mathcal K} (s)}+    A(s)v+ h.o.t \endaligned\end{equation}

\begin{proposition}\label{prop:41} Let $\alpha:\mathbb M\to \mathbb R^3$ be an immersion of class $ C^r $, $ r\geq   6$ and $c$ be a maximal (resp. minimal)
 geometric mean curvature cycle of $\alpha$,
 parametrized by arc length $s$ and  of total length $L$.
 Then the derivative of the Poincar\'e map $\pi_\alpha$ associated to $c$ is given by:

$$ln\pi_{\alpha}^\prime(0)=  \int_0^L\left[\frac{[\sqrt{{\mathcal K}}]_v  }{2\tau_g} +\frac{ k_g}{\tau_g}({\mathcal H}- \sqrt{{\mathcal K}})   \right]ds. \eqno$$

\noindent Here $\tau_g=\sqrt[4]{\mathcal K}\sqrt{2{\mathcal H}-2\sqrt{{\mathcal K}}}$ (resp. $\tau_g=-\sqrt[4]{\mathcal K}\sqrt{2{\mathcal H}-2\sqrt{{\mathcal K}}}$ ).
\end{proposition}

\begin{proof} The Poincar\'e map associated to $c$ is the map $\pi_\alpha:\Sigma \to \Sigma$ defined in a transversal section to $c$ such that $\pi_{\alpha}(p)=p$ for $p\in c\cap\Sigma$ and $\pi_{\alpha}(q)$ is the first return of the geometric mean curvature line through $q$ to the section $\Sigma$, choosing a positive orientation for $c$. It is a local diffeomorphism and is defined, in the local chart $(s,v)$ introduced in Lemma \ref{lm:42}, by $\pi_\alpha:\{s=0\}\to \{s=L\}$, $\pi_\alpha(v_0)=v(L,v_0)$, where $v(s,v_0)$ is the solution of the Cauchy problem
$$(g-\sqrt{{\mathcal K}}G)dv^2+2(f-\sqrt{{\mathcal K}}F)dsdv+(e-\sqrt{{\mathcal K}}E)ds^2=0, \quad v(0,v_0)=v_0.$$

Direct calculation gives
that the derivative  of the Poincar\'e map satisfies
 the following linear differential equation:

$$\frac{d}{ds}(\frac{dv}{dv_0})
 = -\frac{N_v}M(\frac{dv}{dv_0})=-\frac{ [ e-\sqrt{{\mathcal K}}E] _v}{  2[ f-\sqrt{{\mathcal K}}F]  }
  (\frac{dv}{dv_0})\eqno$$

Therefore, using equation \ref{eq:1f2f} it results that
$$\frac{ [ e-\sqrt{{\mathcal K}}E] _v}{  2[ f-\sqrt{{\mathcal K}}F]  } = -\frac{\tau_g^\prime}{2\tau_g}-\frac{[\sqrt{{\mathcal K}}]_v  }{2\tau_g} -\frac{k_g}{\tau_g}({\mathcal H}- \sqrt{{\mathcal K}}). $$

Integrating  the equation above  along an arc
$[s_0,s_1]$ of geometric mean
curvature line,  it follows that:

 \begin{equation}\label{eq:47}
\frac{dv}{dv_0}|_{v_0=0}= \frac
{(\tau_g(s_1))^{\frac{-1}2}}{(\tau_g(s_0))^{\frac{-1}2}}exp[
\int_{s_0}^{s_1} \left[\frac{[\sqrt{{\mathcal K}}]_v  }{2\tau_g} +\frac{ k_g}{\tau_g}({\mathcal H}- \sqrt{{\mathcal K}})   \right]ds.\end{equation}

Applying \ref{eq:47} along the geometric mean curvature cycle of length $ L$, obtain
 $$\frac{dv}{dv_0}|_{v_0=0}=  exp[
  \int_0^L\left[\frac{[\sqrt{{\mathcal K}}]_v  }{2\tau_g} +\frac{ k_g}{\tau_g}({\mathcal H}- \sqrt{{\mathcal K}})   \right]ds. $$

From the equation ${\mathcal K}=(eg-f^2)/(EG-F^2)$ evaluated
at $v=0$  it follows that ${\mathcal K}=\sqrt{{\mathcal K}}[2{\mathcal H}-\sqrt{{\mathcal K}}]-\tau_g^2.$
Solving this equation it follows that $\tau_g=\sqrt[4]{\mathcal K}\sqrt{2{\mathcal H}-2\sqrt{{\mathcal K}}}$.
 \noindent This ends the proof. \end{proof}

  \begin{remark}\label{rm:41}  At this point we show how to extend
 the expression for the derivative of the hyperbolicity of
geometric mean curvature cycles established for class  $C^6$ to  class $C^3$
 (in fact we need only class $C^4$).

The expression \ref{eq:47}
 is the derivative of the transition map
for a geometric mean curvature foliation (which at this point is only of class
$C^1$), along an arc of geometric mean curvature line. In fact, this follows by
approximating  the $C^3$ immersion by one of class  $C^6$. The
corresponding
transition map (now of class $C^4$) whose derivative is given by
expression \ref{eq:47}
converges to the original one (in
class $C^1$) whose expression must given by the same integral, since
the functions involved there are the uniform limits of the corresponding
ones for the approximating immersion.
\end{remark}

  \begin{remark}\label{rm:42}  The expression for the derivative of the Poincar\'e map is obtained by the integration
 of a one form along the geometric mean curvature line $\gamma$. In the case of the arithmetic mean curvature   cycles the correspondent
   expression for the derivative is given by:
$\ln\pi^\prime(0)=\frac 12 \int_0^L\frac{{\mathcal H}_v}{\sqrt{{\mathcal
H}^2-{\mathcal K}}}ds.$
This was proved in \cite{m}.

\end{remark}

\begin{proposition}\label{prop:42} Let $\alpha:\mathbb M\to \mathbb R^3$ be an immersion of class $ C^r $, $ r\geq   6$,  and  $c$ be a
maximal  geometric mean curvature cycle of $\alpha$, parametrized by arc length
and of length $L$. Consider a chart $(s,v)$ as in lemma \ref{lm:42} and consider the deformation
$$\beta_\ep(s,v)=\beta(\ep,s,v)=\alpha(s,v)+\ep [\frac{A_1(s)}6v^3]\delta (v)N(s)\eqno$$
 \noindent where $\delta =1$ in neighborhood of $v=0$, with small support and $A_1(s)=\tau_g(s)>0$.

Then $c$ is a geometric mean curvature cycle of $\beta_\ep$ for all $\ep $ small and $c$ is a hyperbolic
geometric mean curvature cycle for $\beta_\ep$, $\ep\ne 0$.
\end{proposition}

\begin{proof}  In the chart $(s,v)$, for the immersion $\beta$, it is obtained that:

$$\aligned
E(s,v) &= 1-2k_g(s)v+h.o.t\\
F(s,v) &=  0+ 0.v+h.o.t\\
G(s,v) &= 1 +0.v+h.o.t\\
e(s,v)&=\sqrt{{\mathcal K} (s)}+v[\tau_g^\prime(s)-2k_g(s){\mathcal H}(s)\; )]+ h.o.t\\
f(s,v) &= \tau_g(s)+ [2{\mathcal H}(s)-\sqrt{{\mathcal K}(s)}]^\prime v+ h.o.t\\
g(s,v) &=2{\mathcal H}(s)-\sqrt{{\mathcal K} (s)}+v[A(s)+\epsilon A_1(s)]+ h.o.t \endaligned\eqno$$

In the expressions above $E=<\beta_s,\beta_s>$, $F=<\beta_s,\beta_v>$, $G=<\beta_v,\beta_v>$,
 \noindent $\; e=<\beta_{ss},N>$,
 $\; f=<N,\beta_{sv}>,\;$  $\; g=<N,\beta_{vv}>$ and $N=\beta_s\wedge \beta_v /\mid\beta_s\wedge\beta_v\mid.$

Therefore  $c$ is a maximal  geometric mean curvature cycle for all $\beta_\ep$ and at $v=0$ it follows from equation $ {\mathcal K} =(eg-f^2)/(EG-F^2)$ that $${\mathcal K}_v=  \ep \sqrt{{\mathcal K}}A_1(s)+ f(k_g, \tau_g, \sqrt{{\mathcal K}},  {\mathcal H})(s).$$

Therefore, assuming $A_1(s)=\tau_g(s)>0$,  it results that, $$\frac{d}{d\ep}(ln\pi^\prime(0))|_{\ep=0}= -\int_0^L
\frac{d}{d\ep}\left(\frac{(\sqrt{{\mathcal K}})_v}{2\tau_g}\right)ds=-\frac 14 L  < 0.  $$
\end{proof}

As a synthesis of propositions \ref{prop:41} and \ref{prop:42},  the following theorem is obtained.

\begin{theorem}\label{th:41}
An immersion $\alpha\in {\mathcal M}^{r,s}({\mathbb M^2)}$, $ r\geq   6$, is  $C^6-$local geometric mean curvature structurally stable at a geometric mean curvature cycle $c$
if only if,
$$  \int_0^L\left[\frac{[\sqrt{{\mathcal K}}]_v  }{2\tau_g} +\frac{ k_g}{\tau_g}({\mathcal H}- \sqrt{{\mathcal K}})   \right]ds \neq 0. $$
\end{theorem}

\begin{proof} Using propositions  \ref{prop:41} and \ref{prop:42}, the  local topological character of the foliation can be changed by small perturbation of the immersion, when the cycle is not hyperbolic.
\end{proof}

\section{ Geometric Mean Curvature Lines near the Parabolic Line}

 Let  $0$ be a parabolic  point of a $ C^r , \; r\geq   6$, immersion  $\alpha$ parametrized in a Monge chart $(x,y)$ by $\alpha(x,y)=(x,y,h(x,y))$, where

\begin{equation}\aligned h(x,y) =&\frac k2 y^2 + \frac a6 x^3+\frac b2 xy^2+\frac d2 x^2y+\frac c6 y^3\\ +&\frac A{24}x^4+\frac B{6} x^3y+\frac{C}4 x^2y^2+\frac D6 xy^3+\frac E{24}y^4+O(5)\endaligned
\end{equation}

The coefficients of the first and second fundamental forms are given by:

\begin{equation}\label{eq:1f2fp}\aligned E(x,y)=& 1+O(4)\\
F(x,y)=& dk xy^2+\frac{bk}2y^3+O(4)\\
G(x,y)=&1+k^2 y^2+2kb xy^2+kcy^3+O(4)\\
e(x,y)=& ax+dy+\frac A2x^2+Bxy+\frac C2y^2-\frac 12 dk^2y^3+O(4)\\
f(x,y)=& dx+ by+\frac B2x^2+Cxy+\frac D2 y^2-\frac 12 dk^2xy^2-\frac 12 bk^2y^3+O(4)\\
g(x,y)=& k+bx+cy+\frac C2x^2+Dxy+\frac 12{( E-k^3)} y^2\\
-&\frac 12 k^2d x^2y-\frac 32 bk^2 xy^2+dk^2xy^2+(\frac b2-c)k^2y^3+O(4)\endaligned
\end{equation}

The Gaussian curvature is given by

\begin{equation}\label{eq:gc}
 \aligned
{\mathcal K}(x,y)=& k(ax+dy)+ \frac 12(Ak+2ab-2d^2)x^2+(Bk+ac-bd)xy\\
+&\frac 12(Ck+2cd-2b^2)y^2+O(3).\endaligned \end{equation}

The coefficients of the quartic differential equation \ref{eq:mgc4} are given by

\begin{equation}\label{eq:ai}
\aligned A_{40}=&-k(ax+dy)+\frac 12(2a^2+2d^2-2ab-Ak)x^2\\
+ & (2ad-Bk+bd-ac)xy
+\frac 12(2b^2+2d^2-2cd-Ck)y^2\\
+&\frac 16(6dB-3Ab)x^3
+ \frac 12(2dA+3dC-Ac)x^2y\\
+&\frac 12(4dB+3Cb-2Bc)xy^2 \\
+&\frac 16(12dk^3+6dC+6bD-3cC-3dE)y^3+O(4)\\
 \\
A_{31} =& 4ad x^2+4(ab+d^2)xy+4bdy^2+ 2Adx^3\\
+&(6dB+2Ab)x^2y+(4bB+6dC)xy^2+2(dD+bC)y^3+O(4)\\
\\
A_{22}=& 6d^2x^2+12bdxy+6b^2y^2+ (3AD+6dB)x^3\\
+&6(BD+bB+4dC)x^2y+(4dD+3CD+12bC)xy^2+6bDy^3+O(4)\\
\\
A_{13}=& 4k(dx+by)+ (2Bk+4bd)x^2+4(Ck+cd+b^2)xy+(2kD+4bc)y^2\\
+& (6BD+2bB+2dC)x^3+(12CD+2Bc+6Cb)x^2y \\
+& (6D^2+4cC+2dE+2bD-4dk^3)xy^2+(2bE+2cD-4bk^3)y^3+O(4)\\
\\
A_{04}=& k^2+k(2b-a)x+k(2c-d)y+\frac 12[-2ab+2b^2+(2C-A)k+2d^2]x^2\\
+&[(2D-B)k+c(2b-a)+bd]xy+
[c^2+\frac 12(2E-C)k-k^4-cd+b^2]y^2\\
+&\frac 16(18CD+6dB-3Ab+6Cb)x^3+\frac 12(3dC-Ac+2cC)x^2y\\
+&\frac 12(2bE+6DE-6Dk^3-2bk^3-2Bc+3Cb)xy^2\\
+&\frac 16(6cE-6ck^3-3cC-3dE+6bD)y^3 + O(4)\endaligned
\end{equation}

\begin{lemma}\label{lm:p1} Let $0$ be a parabolic point and consider the parametrization $(x,y,h(x,y))$ as  above. If $k>0$ and $a^2 + d^2 \ne 0$ then the set of parabolic points is locally a regular curve normal to the vector $(a,d)$ at $0$.

 If $a\ne 0$ the parabolic curve is transversal to the minimal principal direction $(1,0)$.

If $a=0$ then the parabolic curve is tangent to the principal direction given by $(1,0)$ and has quadratic contact with the corresponding minimal principal curvature line if $dk(Ak-3d^2)\ne 0$.

\end{lemma}

\begin{proof} If $a \ne 0$, from the expression of $\mathcal K$ given by equation \ref{eq:gc} it follows that the parabolic line is given by $x=-\frac{d}{a}y+O_1(2)$ and so is transversal to the principal direction $(1,0)$ at $(0,0)$.

 If $a=0$, from the expression of $\mathcal K$ given by equation \ref{eq:gc} it follows that the parabolic line is given by $y=\frac{ 2d^2-Ak}{2dk}x^2+O_2(3)$  and that $y=-\frac{d}{2k}x^2+O_3(3)$ is the principal line tangent to the principal direction $(1,0)$.
Now the condition of quadratic contact $\frac{ 2d^2-Ak}{2dk}\ne -\frac{d}{2k}$ is equivalent to $dk(Ak-3d^2)\ne 0$.\end{proof}

\begin{proposition}\label{prop:oc-p} Let $0$ be a parabolic point and the Monge chart $(x,y)$ as above.

If $a\ne 0$ then the mean geometric curvature lines are transversal to the parabolic curve and the mean curvatures lines are shown in the picture below, the cuspidal case.

If $a=0$ and $\sigma=dk(Ak-3d^2)\ne 0$ then the  mean geometric
curvature lines are shown in the picture below. In fact, if
$\sigma >0$ then the mean geometric curvature lines are folded
saddles. Otherwise,  if $\sigma <0$ then the mean geometric
curvature lines are folded nodes. The two separatrices of these
tangential singularities, folded saddle and folded node, as
illustrate in the Figure 3 below, are called parabolic
separatrices.
\end{proposition}

\begin{figure}[htbp] \label{fig:dgeo}
\begin{center}
\includegraphics[angle=0, width=11cm]{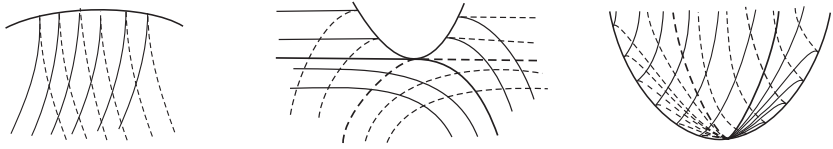}
       \caption {Geometric mean curvature lines near a parabolic point (cuspidal, folded saddle and folded node) and their separatrices}
   \end{center}
 \end{figure}

\begin{proof}
Consider the quartic differential equation
$$H(x,y,p)=A_{04}p^4+A_{13}p^3+A_{22}p^2+A_{31}p+A_{40} =0$$

\noindent where $p=[dx:dy]$ and the Lie-Cartan line field of class $C ^{r-3}$ defined by

$$\aligned x^\prime =&H_p\\
y^\prime =& pH_p\\
p^\prime =&-(H_x+pH_y), \hskip 1cm p=\frac{dy}{dx}\endaligned$$
\noindent where $A_{ij}$ are given by equation \ref{eq:ai}.

If $a\ne 0$ the vector $Y$ is regular and therefore  the mean geometric curvature lines are transversal to the parabolic line and at parabolic points these lines are  tangent to the  principal direction $(1,0)$.

If $a=0$, direct calculation  gives $H(0)=0,\;\; H_x(0)=0, \;\; H_y(0)=-kd, \;\; H_p(0)=0.$

Therefore, solving the equation $H(x,y(x,p),p)=0$ near $0$  it follows,  by the Implicit Function Theorem, that:
$$y=y(x,p)= \frac{2d^2-Ak}{2kd} x^2- \frac{   2Abdk-BAk^2-2d^3b }{2k^2d^2}x^3+O(4).$$

Therefore the vector field
$Y$ given by  the differential equation below
$$\aligned x^\prime =&H_p(x,y(x,p),p)\\
p^\prime =&-(H_x+pH_y)(x,y(x,p),p)\endaligned$$
\noindent is given by

$$\aligned x^\prime =&  4\frac{d^3}k x^3+12d^2x^2p+12dkxp^2+4k^2p^3+O(4)\\
p^\prime =& (Ak-2d^2)x+dkp+O(2).\endaligned
$$

The singular point $0$ is  isolated and the eigenvalues of the
linear part of $Y$ are given by $\lambda_1=0$ and $\lambda_2=dk$.
The correspondent eigenvectors are given by $f_1=(1,(2d^2-Ak)/dk)$
and $f_2=(0,1)$.

Performing the calculations, restricting $Y$ to the center manifold $W^c$ of class $C ^{r-3}$, $T_0W^c=f_1$,  it follows that
$$Y_{c}= -4\frac{(Ak-3d^2)^3}{kd^3}x^3 +0(4)$$

It follows that $0$ is a topological saddle or node  of cubic type provided $\sigma (Ak-3d^2)kd\ne 0$. If $\sigma >0$ then the mean geometric curvature lines are folded saddles and   if $\sigma <0$ then the   mean geometric curvature lines are folded nodes. In the case $\sigma>0$, the center manifold $W^c$ is unique, \cite{sac}, cap. $V$, page 319, and so the saddle separatrices      are well defined. See Figure 4 below.

\begin{figure}[htbp] \label{fig:selano}
\begin{center}
\includegraphics[angle=0, width=10cm]{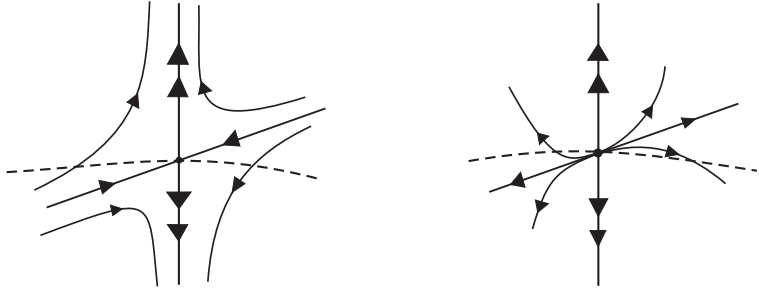}
  \caption {Phase portrait of the vector field $Y$ near singularities }
   \end{center}
 \end{figure}

The reader may find a more complete study of this structure, which can be expressed in the context of normal hyperbolicity,  in the paper of Palis and Takens \cite {pt}.

Notice that due to the constrains of the problem treated here, the non
  hyperbolic saddles and nodes, which in the standard theory would
  bifurcate into three singularities, are actually structurally stable (do
  not bifurcate).

 For a deep analysis of the lost of the hyperbolicity condition and the
  consequent bifurcations,  the reader is addressed to the book of Roussarie
 \cite{rr}.

 \end{proof}

\begin{theorem}\label{th:51}
An immersion $\alpha\in {\mathcal M}^{r,s}({\mathbb M^2)}$, $ r\geq   6$, is  $C^6-$local geometric mean curvature structurally stable at a tangential parabolic point $p$
if only if, the condition $\sigma \neq 0$ in  proposition \ref{prop:oc-p}  holds.

\end{theorem}

\begin{proof}
Direct from Lemma   \ref{lm:p1} and proposition \ref{prop:oc-p}, the  local topological character of the foliation can be changed by small perturbation of the immersion when $\sigma = 0$.
\end{proof}

\section{Examples of Geometric Mean Curvature Configurations}

As mentioned in the Introduction,  no examples of geometric mean curvature foliations are given in the literature, in contrast with the principal and asymptotic foliation. In this section are studied the geometric mean curvature configurations in two classical surfaces: The Torus and the Ellipsoid.
In contrast with the principal case \cite{Sp, St} (but in concordance with the arithmetic mean curvature one \cite{m}) non-trivial recurrence can occur here.

\begin{proposition}\label{prop:mgtoro}
Consider a torus of revolution $T(r,R)$ obtained by rotating a
circle of radius $r$ around a line in the same plane and at a
distance $R$, $R>r$, from its center.
 Define the function
$$\rho=\rho(\frac rR)= 2(\frac rR)^{\frac 34}\int_{-\frac{\pi}2}^{\frac{\pi}2}\frac{ds}{ \sqrt[4]{\cos s}(1+\frac rR\cos s)^\frac 34}.$$
Then the   geometric mean curvature lines on $T(r,R)$, defined in
the elliptic region are all closed or all recurrent according  to
$\rho\in \mathbb Q$ or $\rho \in \mathbb R\setminus\Bbb Q$. Furthermore, both cases occur for appropriate $(r,R)$.
\end{proposition}

\begin{proof} The torus of revolution $T(r,R)$ is parametrized by

 $$\alpha(s,\theta)=((R+r\cos s )\cos \theta ,(R+r\cos s
)\sin \theta, \; r\sin s ).$$

 \noindent Direct calculation shows that $E=r^2,\;$ $\;F=0$, $\;G=[R+r\cos s ]^2,\;$  $e=-r$, $\;f=0$ and
 $\;g=-\cos s(R+r\cos s)$.
  Clearly  $(s,\theta)$ is a principal chart.

 The differential equation of the geometric mean curvature lines, in the principal chart $(s,\theta)$, is given by
$\;\sqrt{eE} ds^2-\sqrt{gG} d\theta^2=0$. This is equivalent to
$$(\frac{ds}{d\theta})^2=\sqrt{\frac{\cos s(R+r\cos s)^3}{r^3}}$$

Solving the equation above it  follows that,

$$\int_{\theta_0}^{\theta_1} d\theta=
  (\frac rR)^{\frac 34}\int_{-\frac{\pi}2}^{\frac{\pi}2}\frac{ds}{ \sqrt[4]{\cos s}(1+\frac rR\cos s)^\frac
34}.$$ So the  two Poincar\'e maps,
$\pi_\pm:\{s=-\frac{\pi}2\}\to
\{s=\frac{\pi}2\}$, defined by $\pi_\pm(\theta_0)=\theta_0 \pm
2\pi\rho(\frac rR)$ have rotation number equal to $\pm \rho(\frac rR)$.
Direct calculations gives that $\rho(0)>0$ and $\rho^\prime(0)<0$.
Therefore,  both the rational and irrational cases occur. This
ends the proof.
\end{proof}

\begin{proposition}\label{prop:mgelipsoide}  Consider an ellipsoid $\mathbb E_{a,b,c}$
with three axes $a>b>c>0$. Then $\mathbb E_{a,b,c}$ have four
umbilic points located in the plane of symmetry orthogonal to
middle axis; they are of the type $G_1$ for geometric mean
curvature lines and of type $D_1$ for the principal curvature
lines.\end{proposition}

\begin{proof} This follows from proposition \ref{prop:2}
and the fact that the arithmetic mean curvature lines have
this configuration, as established in \cite{m}.
\end{proof}

\begin{proposition} Consider an  ellipsoid $\mathbb E_{a,b,c}
$ with three axes $a>b>c>0$. On the ellipse $\Sigma\subset \mathbb
E_{a,b,c}$, containing the four  umbilic points,$\;p_i$,
$i=1,\cdots,4,\;$ counterclockwise oriented, denote by $S_1$
(resp. $S_2$)  the   distance between the adjacent umbilic points
$p_1$ and $p_4$ ( resp. $p_1$ and $p_2$). Define
$\rho=\frac{S_2}{S_1}$.

Then if $\rho\in\Bbb R\setminus \mathbb Q$ (resp. $\rho\in \mathbb
Q$)
 all the  geometric mean curvature lines are recurrent ( resp. all, with the exception of the geometric mean
curvature umbilic separatrices, are closed). See Figure 5 below.
\end{proposition}

\begin{figure}[htbp] \label{fig:selano}
\begin{center}
\includegraphics[angle=0, width=11cm]{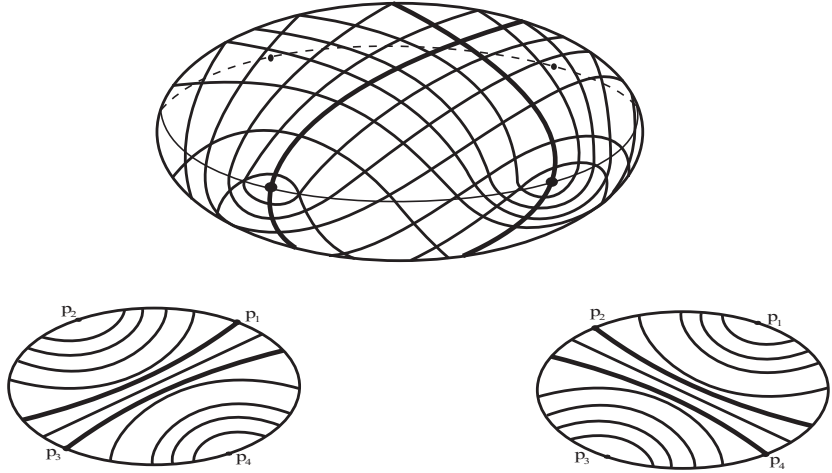}
   \caption {Upper view of geometric mean curvature lines on the
ellipsoid }
   \end{center}
 \end{figure}

\begin{proof}
 The ellipsoid $\mathbb E_{a,b,c}$ belongs to  the triple orthogonal system of surfaces defined by the
one parameter family of quadrics,
$\frac{x^2}{a^2+\lambda}+\frac{y^2}{b^2+
\lambda}+\frac{z^2}{c^2+\lambda}=1$ with $a>b>c>0$, see also
\cite{St}  and \cite{Sp}.

The  following parametrization of $\mathbb E_{a,b,c}$.

$$\alpha(u,v) =
\big(\pm  \sqrt{\frac{M(u,v,a)}{W(a,b,c)}}, \pm \sqrt{\frac{M(u,v,b)}{W(b,a,c)}}, \pm \sqrt{\frac{M(u,v,c)}{W(c,a,b)}}\big)$$
  \noindent  where,

\noindent $M(u,v,w)= {w^2(u+w^2)(v+w^2)}$ and $W(a,b,c)=(a^2-b^2)(a^2-c^2)$,
define  the ellipsoidal  coordinates $(u,v)$ on $\mathbb E_{a,b,c}$, where  $u\in
(-b^2,-c^2)$ and $v\in (-a^2,-b^2)$.

The first fundamental form of $\mathbb E_{a,b,c}$  is given by:

$$I=ds^2=Edu^2+Gdv^2=\frac
14\frac{(u-v)u}{h(u)}du^2+
 \frac 14\frac{(v-u)v}{h(v)}dv^2 $$

The second fundamental form is given by

$$II=
edu^2+gdv^2=\frac{abc(u-v)}{4\sqrt{uv}h(u)
}du^2+\frac{abc(v-u)}{4\sqrt{uv}h(v)}dv^2,$$
\noindent where
$h(x)=(x+a^2)(x+b^2)(x+c^2)$.
  The four umbilic points are $(\pm x_0,0,\pm z_0)= (\pm
a\sqrt{\frac{a^2-b^2}{a^2-c^2}},0,\pm
c\sqrt{\frac{c^2-b^2}{c^2-a^2}}\;)$.

The differential equation of the geometric mean curvature lines is
given by:
$$(\frac{dv}{du})^2=\sqrt{\frac{eE}{gG}}=\sqrt{\frac{-\frac{u}{h(u)}}{-\frac{v}{h(v)}}}.$$

Define $d\sigma_1=\sqrt[4]{-\frac{u}{h(u)}}du$ and
$d\sigma_2=\sqrt[4]{-\frac{v}{h(v)}}dv$. By integration, this leads to the chart $(\sigma_1, \sigma_2)$,
in which the
differential equation of the geometric mean curvature lines is
given by
$$d\sigma_1^2-d\sigma_2^2=0.$$

On the ellipse $\Sigma=\{(x,0,z) | \frac{
x^2}{a^2}+\frac{z^2}{c^2}=1\}$ the distance between the umbilic
points $p_1=(x_0,0,z_0)$ and $p_4=(x_0,0,-z_0)$ is given by
$S_1=\int_{-b^2}^{-c^2}\frac{\sqrt[4]{-v}}{\sqrt{h(v)}}dv$ and
that between the umbilic points $p_1=(x_0,0,z_0)$ and
$p_2=(-x_0,0,z_0)$
 is given by
$S_2=\int_{-a^2}^{-b^2}\frac{\sqrt[4]{-u}}{\sqrt{h(u)}  }du$.

It is obvious that  the ellipse $\Sigma$ is the union of four
umbilic points and the   four principal umbilical separatrices for the
principal foliations.  So $\Sigma\backslash\{p_1,p_2,p_3,p_4\}\;$
is a transversal section of both geometric mean curvature
foliations. The differential  equation of the geometric mean
curvature lines in the principal chart $(u,v)$ is given by
$\sqrt{eE}du^2-\sqrt{gG}dv^2=0$,  which is equivalent to $(\sqrt[4]{
eE} du)^2=(\sqrt[4]{gG}dv)^2$, which amounts to $d\sigma_1=\pm
d\sigma_2$. Therefore near the  umbilic point $p_1$ the geometric
mean curvature lines with a geometric mean curvature umbilic
separatrix contained in the region $\{y>0\}$ define a the transition
map $\sigma_+:\Sigma\to \Sigma$ which is an isometry, reversing
the orientation, with $\sigma_+(p_1)=p_1$. This follows because in
the principal chart $(u,v)$ this
  map is defined by
$\sigma_+:\{u=-b^2\}\to \{v=-b^2\}$ which satisfies the
differential equation $\frac{d\sigma_2}{d\sigma_1}=-1$. By
analytic continuation it results that $\sigma_+$  is an orientation reversing isometry,
  with two fixed points $\{p_1,\;p_3\}$. The
geometric reflection $\sigma_-$, defined in the region $y<0$ have
the two umbilics $\{p_2,\;p_4\}$ as fixed points.

So on the ellipse
parametrized by arclength defined by $\sigma_i$,
  the Poincar\'e return map $\pi_1:\Sigma\to \Sigma$ (
composition of two isometries $\sigma_+$ and $\sigma_-$) is a
rotation with rotation number given by $\frac{S_2}{S_1}$.

Analogously for the other geometric mean curvature foliation,
with the Poincar\'e return map given by $\pi_2=\tau_+\circ \tau_-$,
where $\tau_+$ and $\tau_-$ are two isometries having respectively
$\{p_2,p_4\}$ and $\{p_1,p_3\}$ as fixed points. \end{proof}

\section {Geometric Mean Curvature Structural Stability}

In this section the results of sections 3, 4 and 5 are put together to provide sufficient conditions for geometric mean curvature stability, outlined below.

\begin{theorem} \label{th:sta}
The set of immersions ${\mathcal G}_i({\mathbb M^2}), i=1,\; 2$ which satisfy
conditions $i$), ... , $v$) below
 are i-$C^s$-mean curvature structurally stable and ${\mathcal G}_i, i=1,\; 2$ is open in ${\mathcal M}^{r,s}({\mathbb M^2)}, \; r\geq  s\geq 6$.
\begin{itemize}

\item[i)]  The parabolic curve is regular : ${\mathcal K}=0$ implies $d{\mathcal K}   \neq 0$
and the tangential singularities are saddles and nodes.

\item[ii)]  The umbilic points are of type $G_i$, $i=1,\; 2,\;3$.

\item[iii)]  The geometric mean curvature cycles of $\mathbb G_{\alpha,i}$ are hyperbolic.

\item[iv)]  The  geometric mean curvature foliations $\mathbb G_{\alpha,i}$  has no separatrix connections. This means that there is no  geometric mean curvature line joing two umbilic or tangential parabolic singularities and being separatrices at both ends. See  propositions \ref{prop:2} and \ref{prop:oc-p}

\item[v)] The limit set of every leaf of $\mathbb G_{\alpha,i}$ is    a parabolic point, umbilic point or a geometric mean curvature cycle.
\end{itemize}

\end{theorem}

\begin{proof} The openness of ${\mathcal G}_i({\mathbb M^2})$   follows from the local structure of the geometric mean curvature lines near the umbilic points $G_i$, $i=1,2,3$, near the geometric mean curvature cycles and by the absence of umbilic geometric mean curvature separatrix connections and the absence of recurrences.
 The  equivalence can be performed by the method of canonical regions and their continuation as was done in
\cite{gs1, gs2} for principal lines,  and in \cite{a2}, for asymptotic lines. \end{proof}

 Notice that   Theorem \ref{th:sta} can be reformulated
so as to give the mean geometric stability of the configuration rather than
that of the  separate  foliations. To this end it is necessary to consider the folded extended lines, that  is  to consider the line  of one foliation that arrive  at the parabolic set at a given transversal point as continuing through the line  of the other foliation leaving the parabolic set at  this point, in a sort of ``billiard". This gives raise to the extended folded cycles and separatrices that must be preserved by the homeomorphism mapping simultaneously the two foliations.

Therefore the third, fourth and fifth hypotheses above should be modified as follows:
\begin{itemize}
\item[iii')] the extended folded periodic cycles should be hyperbolic,
\item[iv')] the extended folded separatrices should be disjoint,
\item[v')] the limit set of extended lines should be umbilic points, parabolic singularities and extended folded cycles.
\end{itemize}

 The class of immersions which verify the extended five conditions i), ii), iii'), iv'), v') of a compact and oriented manifold $\mathbb M^2$ will be denoted by ${\mathcal G}({\mathbb M^2})$.

This procedure has been adopted by the authors in the case of asymptotic lines by the suspension operation in order to pass from the foliations to the configuration and properly formulate the stability results.  See \cite{a2}.

\begin{remark} In the space   of convex  immersions ${\mathcal M}^{r,s}_c({\mathbb S}^2)$ ( ${\mathcal K}_\alpha>0$),  the sets  ${\mathcal G}({\mathbb S^2})$ and ${\mathcal G}_1({\mathbb S^2})\cap {\mathcal G}_2({\mathbb S^2}) $ coincide. \end{remark}

 The genericity result involving the five conditions above  is formulated now.

\begin{theorem} The sets ${\mathcal G}_i, \; i=1,\; 2$ are dense in ${\mathcal M}^{r,2}({\mathbb M^2)},\; r\geq  6$.

In the space  ${\mathcal M}^{r,2}_c({\mathbb S}^2)$ the set ${\mathcal G}({\mathbb S^2})$ is dense.
\end{theorem}

The main ingredients for the proof of this theorem are the Lifting
and   Stabilization Lemmas, essential for the achievement of
condition five. The conceptual background for this approach goes
back to the works of Peixoto and Pugh.

The elimination of non-trivial recurrences -- the so called  ``Closing Lemma Problem"-- as a step to achieve condition $v)$ is by far the  most difficult of the steps in the proof.
See the book of Palis and
Melo, \cite{pm}, for a presentation of these ideas in the case of
vector fields on surfaces.

 The proof of theorem above will be postponed to a forthcoming paper \cite{mm}. It involves   technical details that are closer to the proofs of genericity theorems given by
Gutierrez and Sotomayor, \cite{gsln, gs2}, for principal curvature lines  and by Garcia and Sotomayor, \cite{m}, for arithmetic mean  curvature lines.

 \vskip 1cm

\author{\noindent Jorge Sotomayor\\Instituto de Matem\'{a}tica e Estat\'{\i}stica,\\Universidade de S\~{a}o Paulo, \\Rua do Mat\~{a}o 1010, Cidade Universit\'{a}ria, \\CEP 05508-090, S\~{a}o Paulo, S.P., Brazil  \\ \\
\\ Ronaldo Garcia\\Instituto de Matem\'{a}tica e Estat\'{\i}stica,\\Universidade Federal de Goi\'as,\\CEP 74001-970, Caixa Postal 131,\\Goi\^ania, GO, Brazil}

\end{document}